\documentclass[12pt]{article}
\usepackage{a4, amssymb, graphicx, exscale}
\usepackage{makeidx, multicol}

\newenvironment{evlist}[2]{
\begin{list}{}{
\setlength{\topsep}{0.5ex plus0.2ex minus0.1ex} 
\setlength{\leftmargin}{#1}
\setlength{\itemsep}{#2 plus0.2ex}
\setlength{\listparindent}{0pt}
\setlength{\parsep}{0ex plus0.2ex} }}
{\end{list}}

\newcommand{\definition}[1]{\textit{#1}}

\newsavebox{\ttt}
\sbox{\ttt}{}
\begin{document}
\title{Coupled one-dimensional dynamical systems}
\author{Chris Preston}
\date{\small{}}
\maketitle

\thispagestyle{empty}

The aim of this note is to bring attention to a simple class of discrete dynamical systems exhibiting some complex behaviour.
Each of these systems is defined as a self-mapping of the unit square $I^2 = I \times I$ (with $I$ the closed unit interval $[0,1]$) 
and is obtained by coupling two families of self-mappings of the interval $I$.  
There is no real mathematics to be found here (in the sense of results stated and proved)
and in fact there is an almost complete lack of precise statements. 
The only thing on offer is the definition of the mappings and a few nice pictures showing examples of their asymptotic limit sets.
I have written a JavaScript program, accessible at
\texttt{www.math.uni-bielefeld.de/\~{}preston/iterates.html},
which can be used to `discover' more about these mappings. The program might prove to be helpful for anyone interested in doing this.

A mapping $h : Z \to Z$ of a set $Z$ into itself will be considered as a discrete dynamical system in the usual way: The set $Z$ is 
the state space and if the system is in state $z$ at time $n$ then $h(z)$ is the state at time $n + 1$. For each $z \in Z$ the sequence 
$\{z_n\}_{n \ge 0}$, where $z_0 = z$ and $z_{n+1} = h(z_n)$ for all $n \ge 0$, then describes the successive states of the system given 
that the system started in state $z$ at time $0$. This sequence is called the \definition{orbit of $z$ under $h$}.

For various classes of mappings it is expected (or hoped) that some kind of limit set should exist. 
Assuming $Z$ is a metric space with metric $d$ this means that a closed subset $S$ of $Z$ should exist such that the orbit 
$\{z_n\}_{n \ge 0}$ of each `typical' point $z$ converges to the whole of $S$. An orbit $\{z_n\}_{n \ge 0}$ 
converging to the whole of $S$ means that the following two statements hold:

\begin{evlist}{22pt}{5pt}
\item[--]
The orbit eventually comes arbitrarily close to the set $S$. More precisely,
for each $\varepsilon > 0$ there exists $m \ge 0$ such that $d(z_n,S) < \varepsilon$ for all $n \ge m$
(i.e., for each $n \ge m$ there exists $y_n \in S$ with $d(z_n,y_n) < \varepsilon$). 
\item[--]
Each $y \in S$ is an accumulation point of the orbit. More precisely, for each $\varepsilon > 0$ there exist infinitely many indexes $n$ 
with $d(y,z_n) < \varepsilon$.
\end{evlist}

(Note that if there is a closed set $S$ satisfying these two conditions then it is uniquely determined by the orbit 
$\{z_n\}_{n \ge 0}$.)
That this should hold for a `typical' point $z$ means that the set of points for which it does not hold should be negligible
in either a measure-theoretical or topological sense (i.e., it should have measure zero with respect to an appropriate measure on $Z$
or be of the first category in the terminology of the Baire category theorem).

Suppose that a limit set $S$ is thought to exist for some explicitly given mapping $h : Z \to Z$. Then an approximation to $S$ can be 
displayed on a monitor as follows: Start with some random initial point $z$ (which it is hoped will be 
`typical') and consider the orbit $\{z_n\}_{n \ge 0}$ of $z$ under $h$.
Compute the first $N$ terms of this orbit, where $N$ is large enough so that $z_n$ is less than one pixel from $S$ for all $n \ge N$
and then display the next $M$ terms of the orbit on the monitor, where $M$ is chosen so that the pixels of
the displayed points more-or-less fill out the set $S$. Of course, these statements concerning the choice of $M$ and $N$ are
extremely vague and in practice it is necessary to resort to trial and error to find suitable values for them.

Conversely, suppose it is not known whether a limit set exists, but when the above procedure is carried out for many different 
initial points $z$ and for ever larger values for $N$ the resulting image is always the same.
Then this can be interpreted as evidence that a limit set exists and that the image is an approximation to the limit set.
In this sense there is strong evidence that each of
the two-dimensional dynamical systems which we now introduce possesses a limit set.

As already stated, the dynamical systems we are interested in are defined in the unit square $I^2$ and are obtained by coupling two 
families of self-mappings of the interval $I$.
Each such family of self-mappings of $I$ will be given by a mapping $f : I^2 \to I$, where for each $p \in I$ the mapping 
$x \mapsto f(p,x)$ is the 
self-mapping of $I$ corresponding to the parameter value $p$. We write $f_p(x)$ instead of $f(p,x)$ to emphasise 
that the first argument is to be regarded as a parameter.
The prototypical example is the \definition{logistic map} (or \definition{family}) $\ell$ with 
\[\ell_p(x) = 4px(1 - x)\;.\]
Another typical example is the \definition{tent map} (or \definition{family}) $t$ with 
\[t_p(x) = p(1 - |2x - 1|) \;.\]

As well as two families of self-mappings of $I$ we also require two couplers. A coupler is just a mapping $c : I \to I$, for example there a 
simple linear coupler with $c(x) = b + rx$, where $b, r$ are constants with $b \ge 0$, $r \ge 0$ and $b + r \le 1$.
This coupler will be referred to as the \definition{linear+ coupler} and the parameters $b$ and $r$ as its \definition{base} and 
\definition{rate}.

Given two families of self-mappings $f,\, g : I^2 \to I$ and two couplers $c,\,d : I \to I$ the coupled system is then the 
mapping $h : I^2  \to I^2$ defined by
\[      h(x,y) = (f_{c(y)}(x), g_{d(x)}(y)) \]
for all $(x,y) \in I^2$. The first component of $h(x,y)$ is thus obtained by applying the family $f$ to the
first argument $x$ with the parameter value given by applying the coupler $c$ to the second argument $y$. In the same way
the second component of $h(x,y)$ is obtained by applying the family $g$ to the
second argument $y$ with the parameter value given by applying the coupler $d$ to the first argument $x$. 

To be more explicit, in the following let us take both $f$ and $g$ to be the logistic map $\ell$ and $c$ and $d$ to be linear+ couplers, 
$c$ with base $b$ and rate $r$ and $d$ with base $b'$ and rate $r'$. Then the mapping $h$ is given by
\[
 h(x,y) = (4(b + ry)x(1 - x), 4(b' + r'x)y(1 - y))
\]
and so $h$ depends on the four parameters $b$, $r$, $b'$ and $r'$.
Carrying out the procedure described above for this mapping $h$ with $b = b' = 0.4$ and $r = r' = 0.6$ results in the following image: 

\medskip

\begin{center}
\includegraphics[width=10cm]{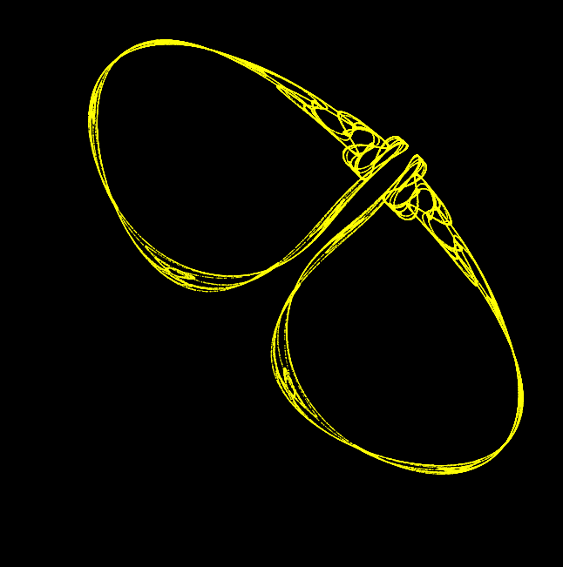}
\end{center}

\medskip

The values $N = 1000000$ and $M = 100000$ were used here and the initial value was $(0.7,0.6)$.
Changing the initial value or increasing $N$ does not alter the image. The value for $M$ seems about right, but this
is rather subjective and also depends on the monitor.  

Now what is more interesting than the limit set for a particular mapping is how the limit set behaves as a function of 
the parameters $b$, $r$, $b'$ and $r'$. If we choose a line segment (or more general curve)
in this four-dimensional parameter space and compute the limit sets along a grid of points on the curve then we obtain
a succession of images which can be made to appear as a `video'.

Consider the simple case in which $b'$ and $r'$ are fixed with $b' = 0.4$ and $r' = 0.6$, where $r = 1 - b$ and
$b$ runs from $0$ to $1$. In other words we are considering the line segment $s \mapsto (s,1-s,0.4,0.6)$ in parameter space.
Here is a selection of limit sets from points lying on this line segment:

\begin{center}
\includegraphics[width=4cm]{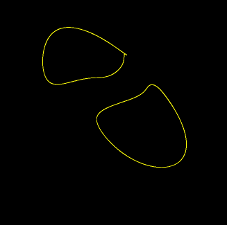}
\includegraphics[width=4cm]{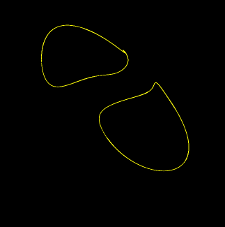}
\includegraphics[width=4cm]{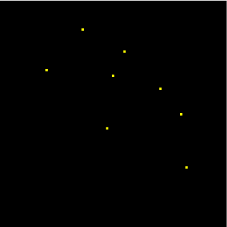}

$\ b=0.00\quad r=1.00$ \qquad $b=0.05\quad r=0.95$ \qquad $b=0.10\quad r=0.90$

\medskip

\includegraphics[width=4cm]{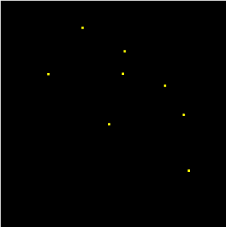}
\includegraphics[width=4cm]{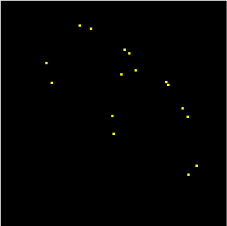}
\includegraphics[width=4cm]{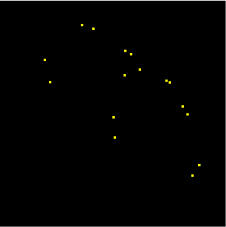}

$\ b=0.15\quad r=0.85$ \qquad $b=0.20\quad r=0.80$ \qquad $b=0.25\quad r=0.75$

\medskip

\includegraphics[width=4cm]{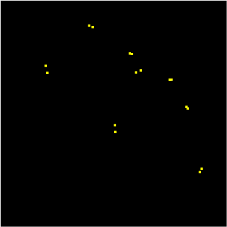}
\includegraphics[width=4cm]{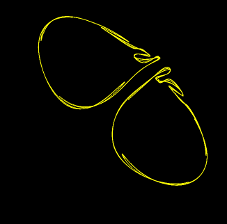}
\includegraphics[width=4cm]{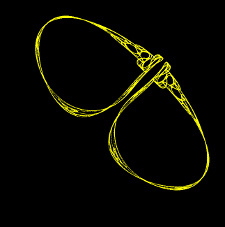}

$\ b=0.30\quad r=070$ \qquad $b=0.35\quad r=0.65$ \qquad $b=0.40\quad r=0.60$

\medskip
\includegraphics[width=4cm]{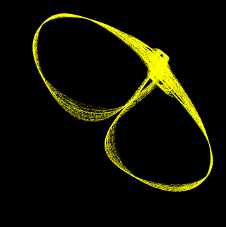}
\includegraphics[width=4cm]{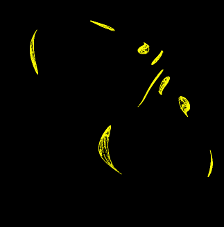}
\includegraphics[width=4cm]{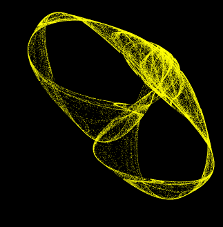}

$\ b=0.45\quad r=0.55$ \qquad $b=0.50\quad r=0.50$ \qquad $b=0.55\quad r=0.45$

\end{center}

\begin{center}

\includegraphics[width=4cm]{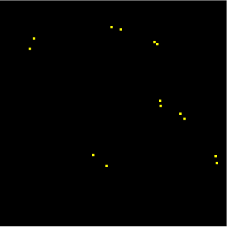}
\includegraphics[width=4cm]{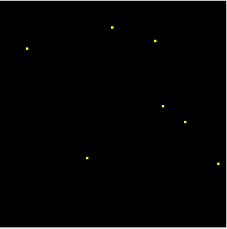}
\includegraphics[width=4cm]{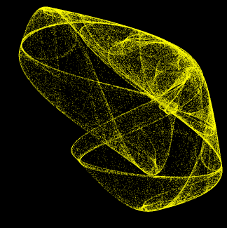}

$\ b=0.60\quad r=1.85$ \qquad $b=0.65\quad r=0.80$ \qquad $b=0.70\quad r=0.75$

\medskip

\includegraphics[width=4cm]{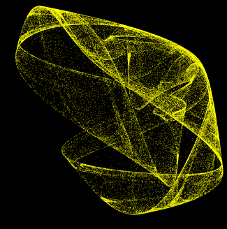}
\includegraphics[width=4cm]{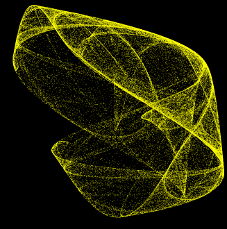}
\includegraphics[width=4cm]{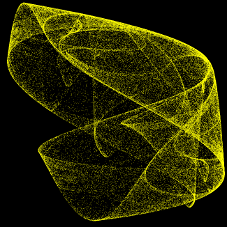}

$\ b=0.75\quad r=0.25$ \qquad $b=0.80\quad r=0.20$ \qquad $b=0.85\quad r=0.15$

\medskip

\includegraphics[width=4cm]{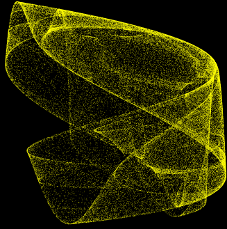}
\includegraphics[width=4cm]{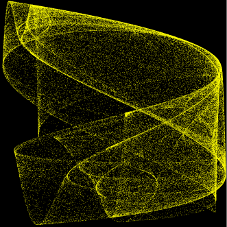}
\includegraphics[width=4cm]{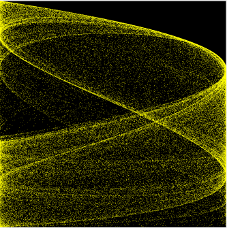}

$\ b=0.90\quad r=0.10$ \qquad $b=0.95\quad r=0.05$ \qquad $b=1.00\quad r=0.00$

\end{center}

\bigskip

The images shown above are taken from a uniform grid of points along the line segment and give a very rough impression
of some of the different kinds of limit sets which can occur. (The images where the limit set is a finite periodic cycle
have been manipulated by making the points larger. Before this was done they were rather hard to see.)
In order to get a more complete picture, however, it is necessary
to consider a much finer grid of say 1000 points (which can be done using the JavaScript program).

The qualitative behaviour of the limit sets along any line segment of the form $s \mapsto (s,1-s,b',r')$ 
with $b' + r' = 1$ seems to be essentially the same as that along the line segment $s \mapsto (s,1-s,0.4,0.6)$.
Moreover, if one or both of the logistic maps is replaced with the tent map and the limit sets are considered along one of
these line segments then the resulting succession of images is basically similar to that from the  original example.

We now look at something slightly different. 
Again consider two families of self-mappings $f,\, g : I^2 \to I$ and two couplers $c,\,d : I \to I$.
Then there is another way of obtaining a coupled system: This uses the mapping $h' : I^2  \to I^2$ defined by
\[ h'(x,y) = (f_{c(y)}(x), g_{d(f_{c(y)}(x))}(y)) \]
instead of using $h$. To see what is going on here let us write the definitions of $h$ and $h'$ one under the other in 
the following form: 
\[h(x,y) = (x',y'),\ \mathrm{where}\ x' = f_{c(y)}(x)\ \mathrm{and}\ y' = g_{d(x)}(y),\]
\[h'(x,y) = (x',y'),\ \mathrm{where}\ x' = f_{c(y)}(x)\ \mathrm{and}\ y' = g_{d(x')}(y).\]
From this it can perhaps be seen how the mapping $h'$ could be `discovered'  as a result of
making a common programming mistake when writing a program for $h$.

The reason for looking at this new class of mappings is that there is again strong evidence that the mappings
possess limit sets, but the behaviour of these limit sets along the usual line segments seems to be very different to
that which was seen above.
As before, let us restrict our attention to the case in which both $f$ and $g$ are the logistic map $\ell$ and $c$ and $d$ are linear+ 
couplers, $c$ with base $b$ and rate $r$ and $d$ with base $b'$ and rate $r'$. Then the mapping $h'$ is given by
\[
 h'(x,y) = (4(b + ry)x(1 - x), 4(b' + 4r'(b + ry)x(1 - x))y(1 - y))\;.
\]
Consider the same line segment $s \mapsto (s,1-s,0.4,0.6)$ in parameter space
as above, i.e., with $b'$ and $r'$ held fixed and taking the values $0.4$ and $0.6$ respectively, with $r = 1 - b$ and
with $b$ running from $0$ to $1$.
Here are the limit sets taken at the same uniform grid of points along the line segment:

\newpage

\begin{center}
\includegraphics[width=4cm]{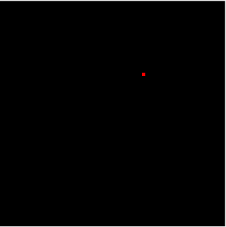}
\includegraphics[width=4cm]{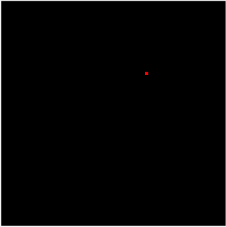}
\includegraphics[width=4cm]{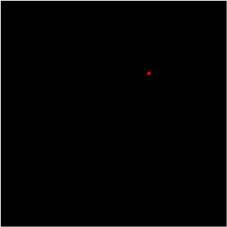}

$\ b=0.00\quad r=1.00$ \qquad $b=0.05\quad r=0.95$ \qquad $b=0.10\quad r=0.90$

\medskip

\includegraphics[width=4cm]{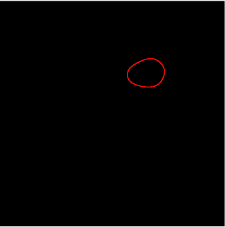}
\includegraphics[width=4cm]{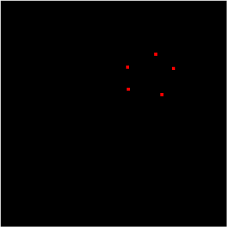}
\includegraphics[width=4cm]{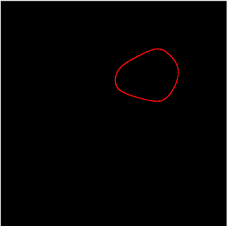}

$\ b=0.15\quad r=0.85$ \qquad $b=0.20\quad r=0.80$ \qquad $b=0.25\quad r=0.75$

\medskip

\includegraphics[width=4cm]{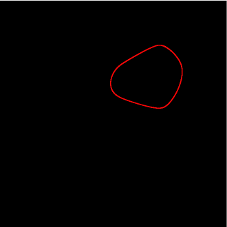}
\includegraphics[width=4cm]{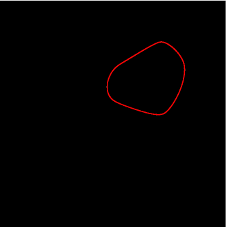}
\includegraphics[width=4cm]{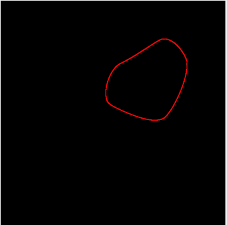}

$\ b=0.30\quad r=070$ \qquad $b=0.35\quad r=0.65$ \qquad $b=0.40\quad r=0.60$

\medskip
\includegraphics[width=4cm]{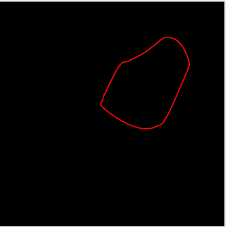}
\includegraphics[width=4cm]{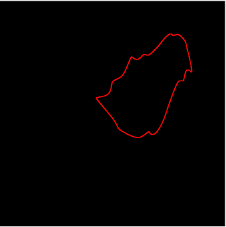}
\includegraphics[width=4cm]{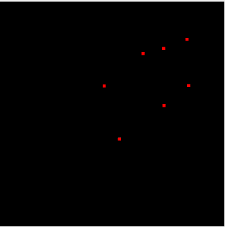}

$\ b=0.45\quad r=0.55$ \qquad $b=0.50\quad r=0.50$ \qquad $b=0.55\quad r=0.45$

\end{center}

\begin{center}

\includegraphics[width=4cm]{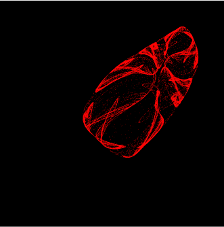}
\includegraphics[width=4cm]{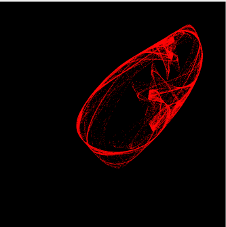}
\includegraphics[width=4cm]{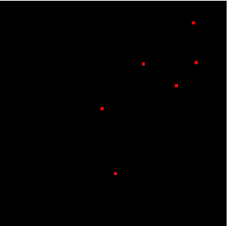}

$\ b=0.60\quad r=1.85$ \qquad $b=0.65\quad r=0.80$ \qquad $b=0.70\quad r=0.75$

\medskip

\includegraphics[width=4cm]{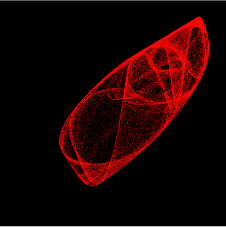}
\includegraphics[width=4cm]{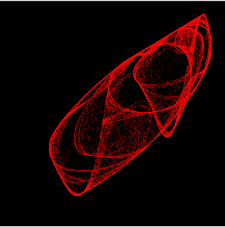}
\includegraphics[width=4cm]{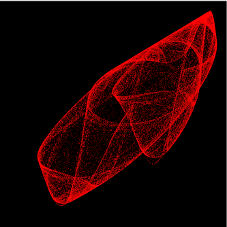}

$\ b=0.75\quad r=0.25$ \qquad $b=0.80\quad r=0.20$ \qquad $b=0.85\quad r=0.15$

\medskip

\includegraphics[width=4cm]{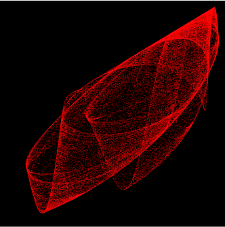}
\includegraphics[width=4cm]{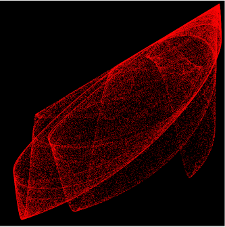}
\includegraphics[width=4cm]{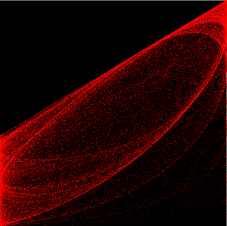}

$\ b=0.90\quad r=0.10$ \qquad $b=0.95\quad r=0.05$ \qquad $b=1.00\quad r=0.00$

\end{center}

\bigskip
(Again, where the limit set is a finite periodic cycle the points have been enlarged to make them easier to see.)
As in the first case, the qualitative behaviour of the limit sets along any of the lines segments considered here
seems to be essentially the same as that along the line segment $s \mapsto (s,1-s,0.4,0.6)$.
Moreover, replacing one or both of the logistic maps with the tent map does not change
the resulting succession of images beyond recognition.

\bigskip

{\sc Fakult\"at f\"ur Mathematik, Universit\"at Bielefeld}\\
{\sc Postfach 100131, 33501 Bielefeld, Germany}\\
\textit{E-mail address:} \texttt{preston@math.uni-bielefeld.de}\\

\end{document}